\documentclass[10pt]{amsart}
\usepackage{amsthm}
\usepackage{amsmath}
\usepackage{amsfonts}
\usepackage{amssymb}

\def\vbar{\overline{v}}

\def\complexes{\mathbb{C}}
\def\reals{\mathbb{R}}
\def\integers{\mathbb{Z}}
\def\rationals{\mathbb{Q}}
\def\U{\mathcal{U}}
\def\H{\dagger}

\newcommand{\beq}{\begin{equation}}
\newcommand{\eeq}{\end{equation}}

\newcommand{\bea}{\begin{eqnarray}}
\newcommand{\eea}{\end{eqnarray}}
\newcommand{\ba}{\begin{array}}
\newcommand{\ea}{\end{array}}
\newcommand{\bean}{\begin{eqnarray*}}
\newcommand{\eean}{\end{eqnarray*}}

\newcommand{\bit}{\begin{itemize}}
\newcommand{\eit}{\end{itemize}}
\newcommand{\ben}{\begin{enumerate}}
\newcommand{\een}{\end{enumerate}}

\begin{document}

\renewcommand{\baselinestretch}{1.1}

\newtheorem{Theorem}{Theorem}[section]
\newtheorem{Lemma}[Theorem]{Lemma}
\newtheorem{Proposition}[Theorem]{Proposition}
\newtheorem{Corollary}[Theorem]{Corollary}

\theoremstyle{remark}
\newtheorem{Remark}[Theorem]{Remark} 
\newtheorem{Notation}[Theorem]{Notation}
\newtheorem{Definition}[Theorem]{Definition}
\newtheorem{Note}[Theorem]{Note}
\newtheorem{Example}[Theorem]{Example}

\author{B.A. Sethuraman}
\address{Department of Mathematics\\California State University
Northridge\\Northridge CA 91330\\U.S.A.}
\title{Division Algebras and Wireless Communication}
\email{al.sethuraman@csun.edu}
\thanks{The author is supported in part by NSF grant DMS-0700904. The author wishes to thank P. Vijay Kumar for innumerable discussions during the preparation of this article: his counsel was invaluable, his patience monumental.}

\maketitle


\section*{}  \label{secn:intro}

The aim of this note is to bring to the attention of a wide
mathematical audience the recent application of division algebras to
wireless communication.  The application occurs in the context of
communication involving multiple transmit and receive antennas, a
context known in engineering as MIMO, short for multiple input, 
multiple output. 
While the use of multiple receive antennas goes
back to the time of Marconi, the basic theoretical framework for
communication using multiple transmit antennas was only published
about ten years ago.  The progress in the field has been quite
rapid, however, and MIMO communication is  widely credited with
being one of the key emerging areas in
telecommunication.  Our focus here will be on one aspect of
this subject: the formatting of transmit information for optimum
reliability.

Recall that a division algebra is an (associative) algebra with a
multiplicative identity in which every nonzero element is
invertible.  The center of a division algebra is the set of
elements in the algebra that commute with every other element in
the algebra; the center is itself just a commutative field, and
the division algebra is naturally a vector space over its center.
We consider only division algebras that are finite-dimensional as
such vector spaces.  Commutative fields are trivial
examples of these division algebras, but they are by no means the
only ones: for instance, class-field theory tells us that over any
algebraic number field $K$, there is a
rich supply of noncommutative division algebras whose center is
$K$ and are finite-dimensional over $K$.

Interest in MIMO communication began with the papers
\cite{TSC,FG,Tel,GFBK} where it was established that MIMO  wireless transmission could be used 
both to decrease the probability of error as well as to increase the amount of
information that can be transmitted.  This caught the attention of telecommunication operators, particularly since MIMO communication does not require additional resources in the form of either a larger slice of the radio spectrum or else increased transmitted power.

The
basic setup is as follows: Complex numbers $r e^{\imath\phi}$, encoded as the amplitude ($r$)
and phase ($\phi$) of a radio wave, are sent from $t$ transmit antennas (one
number from each antenna), and the encoded signals are then received
by $r$ receive antennas. The presence of obstacles in the
environment such as buildings causes attenuation of  the signals; in
addition, the signals are reflected several times  and interfere
with one another.  The combined degradation of the signals is commonly referred to as ``fading'', and achieving reliable communication in the presence of fading has been the most challenging aspect of wireless communication.  The received and transmitted signals are
modeled by the relation
\begin{equation*} \label{eqn:block_model_general_r_t}
Y_{r\times 1} = \theta H_{r\times t} X_{t\times 1} + W_{r\times 1}
\end{equation*}
where $X$ is a $t\times 1$ vector of transmitted signals, $Y$ is an
$r\times 1$ vector of received signals, $W$ is an $r\times 1$ vector
of additive noise,  $H$ is an $r\times t$ matrix that models the
fading, and $\theta$ is a real number chosen to  normalize
the transmitted signals so as to fit the available power. Under the most commonly
adopted model, the entries of the noise vector $W$  and the channel
matrix $H$ are assumed to be Gaussian complex random variables that are
independent and identically distributed  with zero mean. (A Gaussian complex random variable is one of the form $w = x + \imath y$ where $x$ and $y$ are real Gaussian random variables that are independent and have the same mean and variance. The modulus of such a random variable, and in particular the magnitude of
of each fading coefficient $h_{ij}$, is then Rayleigh distributed.
This model is hence also known as the Rayleigh fading channel model.) It
is the presence of fading in the channel that distinguishes this
model from more classical channels, where the primary source of
disturbance is the additive Gaussian noise $W$.

 The
transmission typically occurs in blocks of length $n$: each
antenna transmits $n$ times, and the receiver waits to receive all
$n$ transmission before processing them. 
A common engineering model is to assume that $r=t=n$, 
so the equation 
above is accordingly modified to read
\begin{equation}
 \label{eqn:block_model} Y_{n\times n} =
\theta H_{n\times n} X_{n\times n} + W_{n\times n}.
\end{equation}

Thus, the $i$th column of $Y$, $X$, and $W$ represent (respectively)
the received vectors, the transmitted information, and the additive
noise from the $i$th transmission.  In the model considered here, the fading characteristics of the channel (i.e., the $h_{i,j}$) are
assumed to be known to the receiver, but not to the
transmitter (this is known as \textit{coherent} transmission).  A measure of the power available during a single transmission from all $n$ antennas, i.e., a single use
of the telecommunication channel, is the \textit{signal-to-noise ratio} (SNR) $\rho$.  Recall that the Frobenius norm $||X||_F$ of $X = (x_{i,j})$ equals $\sqrt{\sum_{i,j} |x_{i,j}|^2}$. Since the power required to send a complex number varies as the square of its modulus, the normalization constant $\theta$ must satisfy $\theta^2 ||X||_F^2 \le n \rho$.

A subset $S$ of the nonzero complex numbers known as the
\textit{signal set} is selected as the alphabet (a common situation
is that $S$ is a finite subset of size $q$ of the nonzero Gaussian integers
$\integers[\imath]- \{0\}$), and a $k$-tuple $ (s_1,s_2,\dots,s_{k})$, $s_i
\in S$, comprises the message that the transmitter wishes to convey
to the receiver.  Thus there are $q^k$ messages in all and it is
assumed that each message is equally likely to be transmitted.  A
\textit{space-time code} is then a one-to-one map $X\colon S^{k}\rightarrow
M_n(\complexes)$; we write  
$\mathcal X$ for $X(S^k)$. The transmitted matrix $\theta X_{n\times n}$
in Equation (\ref{eqn:block_model}) is thus drawn from the set
$\theta {\mathcal X}$ as
$(s_1,s_2,\dots,s_{k})$ vary in $S^k$. 
Often ${\mathcal X}$ itself is referred to as
the space time code.   It is typically assumed that the map $X$ is ``linear in $S^k$'', that is, it is the restriction to $S^k$ of a group homomorphism $\langle S\rangle ^k\rightarrow M_n(\complexes)$, where $\langle S \rangle$ is the additive subgroup of $\complexes$ generated by $S$.

Under the information-theoretic framework developed by Shannon in
1948 (\cite{Shan}) and adopted ever since within the telecommunication community,
the amount of information conveyed by a message in this setting is
equal to $\log_2(q^k)$ ``bits''.  Since this amount of information
is conveyed in $n$ transmissions over the MIMO channel, the rate of
information transmission is then given by $\frac{k}{n}\log_2(q)$
bits per channel use.  When $q$ and $n$ are fixed a priori, the
quantity $k$ serves as a measure of information rate.

Reliability of communication is commonly measured by the probability
$P_e$ of incorrectly decoding the transmitted message at the
receiver.    The pairwise error probability $P_e(i,j)$ is the
probability that message $i$ is transmitted and message $j$ is
decoded. Performance analysis
of MIMO communication systems typically focuses on the pairwise
error probability as it is easier to estimate and also because the
error probability $P_e$ can be upper and lower bounded in terms of
the pairwise error probability.

It was shown in \cite{TSC,GFBK} that for a fixed SNR (i.e., power) $\rho$, in order to keep the
pairwise error probability low, the space time code ${\mathcal X}$ must
meet the two criteria below, of which the first is primary:
\begin{enumerate} \label{initial_perf_crit}
\item \label{fullrank} \textit{Rank Criterion:}
For $(s_1,s_2,\dots,s_{k}) \neq (s'_1,s'_2,\dots,s'_{k})$, $X(s_1,s_2,\dots,s_{k})-X(s'_1,s'_2,\dots,s'_{k}))$ must have full
rank $n$, i.e., it must be invertible.
\item \label{coding_gain} \textit{Coding Gain Criterion:}
For $(s_1,s_2,\dots,s_{k}) \neq (s'_1,s'_2,\dots,s'_{k})$, the
modulus of the determinant of difference
$|\det(\theta X(s_1,s_2,\dots,s_{k})-\theta X(s'_1,s'_2,\dots,s'_{k}))|$ must be
as large as possible.
\end{enumerate}
Clearly, the second criterion comes into play only when the first
criterion has been met.  Note that one cannot arbitrarily scale the matrices $X$ to increase the coding gain because the assumption of fixed $\rho$ along with the relation $\theta^2 ||X||_F^2 \le n \rho$ would cause a corresponding decrease in $\theta$.  Note too that the second criterion shows that one cannot increase the quantity $k$ (a proxy for the rate of information) arbitrarily, as this would create a larger set of matrices $\theta X$ all circumscribed to lie within a sphere of radius  $\sqrt{n\rho}$, which would entail that the determinant of their differences would get smaller, thereby increasing the pairwise error probability.

\section*{Satisfying the Rank Criterion} \label{secn:rc}
The earliest space-time code, for two antennas, was given by an
engineer Alamouti (\cite{Ala}): given an arbitrary signal set $S$, he chose $X\colon S^2 \rightarrow
M_2(\complexes)$ to be
\begin{equation} \label{eqn:alamouti_code}
X(s_1,s_2) = \left(
              \begin{array}{cc}
                s_1 & -\overline{s_2} \\
                s_2 & \overline{s_1} \\
              \end{array}
            \right)
\end{equation} (where $\overline{s_i}$ stands for complex conjugation).
It is easy to
see that the rank criterion is immediately met.
Writing $s_1=u_1+\imath u_2$, $s_2=u_3+\imath u_4$, 
each such matrix can be expressed in the form 
$ X(s_1,s_2) =  \sum_{i=1}^4 u_i A_i $.
The $2 \times 2$ complex matrices $A_j$ are such that for any complex $2 \times 2$ channel matrix $H$, the collection of $2 \times 2$ matrices $\{HA_i\}$ is pairwise orthogonal when regarded as vectors in $\reals^8$ by writing out sequentially the real and imaginary parts of each entry of the $\{HA_i\}$.   The expansion above makes it possible to do a least squares estimation of the $u_j$  from the received matrix $Y$,
also considered as a vector in $\reals^8$ as above, by projecting onto the respective matrices 
$HA_j$ (we will consider this in more detail later).  It is this property that makes the Alamouti code so easy to decode, and not surprisingly, the code has since been adopted into the IEEE 802.11n ``Wireless LAN'' standard.  
In applications, the $\{s_1, s_2\}$ are typically drawn from a subset of $\integers[\imath]\times \integers[\imath]$.

Alamouti's code led
to a furious search among engineers and coding theorists for
generalizations for higher number of antennas.   Much of the early work (see \cite{TJC} for example) focused on combinatorial methods.  
The matrix $X$ in Equation (\ref{eqn:alamouti_code}) is almost unitary: it satisfies 
$X X^{\H} = (s_1\overline{s_1} + s_2\overline{s_2}) I_2$, where the superscript $\H$ stands 
for transpose conjugate, and $I_2$ stands for the $2\times 2$ identity matrix.  
Not surprisingly, early workers (see \cite{TJC} for example) 
sought $n\times n$ matrices $X(s_1,\dots, s_k)$ whose entries come 
from the set $\{\pm s_j, \pm \overline{s_j}, \pm \imath s_j, \pm \imath\overline{s_j},\ j = 1,\dots, k\}$  and 
satisfy 
\begin{equation} \label{eqn:orth_des}
X X^{\H} = (s_1\overline{s_1} +\cdots+ s_k\overline{s_k}) I_n
\end{equation}  
This quickly leads to a necessary condition: the existence of $2k-1$  complex $n\times n$ matrices 
$A_i$ satisfying $A_i^{\H} A_i = I_n$, $A_i^{\H} = -A_i$, and 
$A_i A_j = - A_j A_i$ for $1 \le i < j \le 2k-1$.  These are of course the Hurwitz-Radon-Eckmann matrices, and
 classical results of Hurwitz-Radon-Eckmann (see \cite{Ecm} for instance) severely limits the values of $k$ for which 
such matrices can exists.
If $n=2^a(2b+1)$ then the  Hurwitz-Radon-Eckmann result says that the maximum possible value of $k$ equals $(a+1)$.   Thus $k=n$ if and only if   $n=2$,   $k \le \frac{3n}{4}$ for $n >2$, and $k \le \frac{n}{2}$ for $n>4$. 
  It follows that these generalizations of the Alamouti code transmit too few information symbols for more than two transmit antennas.
 (A similar analysis of the matrices $A_i$ using representation of
  Clifford Algebras was made by Tirkkonen and Hottinenin
  \cite{TH}.)

In 2001, Sundar Rajan, a professor of communication engineering at
the Indian Institute of Science introduced the problem of designing
matrices $X(s_1,\dots,s_k)$ satisfying the rank criterion to this author.
Given his algebraic background, this author could recognize  easily
that matrices arising from embeddings of fields and division
algebras can be utilized to solve this problem.
Let
$f\colon D \rightarrow M_n(\complexes)$ be an embedding, i.e., an
(injective) ring homormorphism of a division algebra $D$
into the $n\times n$ matrices over $\complexes$.  Then for $X_1 = f(d_1)$ and $X_2 =
f(d_2)$ ($X_1 \neq X_2$), $X_1-X_2$ must necessarily be
invertible.  This is because $d_1-d_2$, being a nonzero element of the division algebra $D$, is
automatically invertible, and since $f$ is a homomorphism, the same must also be true of $X_1-X_2$.  
Thus, the matrices in $f(D)$ automatically satisfy the rank criterion.
  Using this observation,  Sundar Rajan, his Ph.D. student Shashidhar,
and this author (\cite{SRS}) proposed several schemes for
constructing space-time codes from various signal sets. For each
signal set $S$ and for each $n$, they constructed suitable division algebras
$D$, suitable embeddings $f\colon D \rightarrow M_n(\complexes)$, and suitable
injective maps $X \colon S^{k} \rightarrow f(D)$, for suitable $k$.

For simplicity of construction in the noncommutative case, the authors of (\cite{SRS})
used \textit{cyclic division algebras} for their codes.  A cyclic
division algebra is constructed from two data: a field extension $K/F$ of
degree $n$ that is Galois with cyclic Galois group $\langle
\sigma\rangle$, and a nonzero element $\gamma \in F$ that
satisfies the property that for any $i=1,\dots,n-1$, $\gamma^i$ is
not a norm\footnote{this is a sufficient
condition to obtain a cyclic division algebra} from $K$ to $F$.  As a $K$-vector space, the algebra is expressible as
$$
D = \bigoplus_{i=0}^{n-1} K u^i
$$
where $u$ is a symbol.  The multiplication in this algebra is
given by the relations $u k = \sigma(k) u$ for all $k\in K$, and
$u^n = \gamma$.  The bilinearity of multiplication along with
these relations then allows us to determine the product of any two
elements of $D$. One can prove that this construction 
indeed yields a division algebra with center $F$. (Such a division algebra is said to be of \textit{index $n$}.)

There is a well-known embedding of such a $D$
into $M_n(K)$ that sends $k_0 + k_1 u + \cdots k_{n-1}u^{n-1}$ to
%
%
\begin{equation}
\left[
\begin{array}{ccccc}\label{cycmat}
 k_0 & \gamma \sigma(k_{n-1})  & \gamma \sigma^2(k_{n-2}) & \dots& \gamma \sigma^{n-1}(k_1)
\\
 k_1 & \sigma(k_0) &\gamma \sigma^2(k_{n-1})& \dots  
&\gamma \sigma^{n-1}(k_2) \\
 k_2 & \sigma(k_1) &\sigma^2(k_0)& \dots &\gamma
\sigma^{n-1}(k_3) \\
 k_3 & \sigma(k_2) &\sigma^2(k_1)& \dots &\gamma
\sigma^{n-1}(k_4)\\
 \vdots & \vdots &\vdots& \vdots&  \vdots\\
  k_{n-2} & \sigma(k_{n-3}) & \sigma^2(k_{n-4})&\dots &
\gamma\sigma^{n-1}(k_{n-1})\\
 k_{n-1} & \sigma(k_{n-2}) & \sigma^2(k_{n-3})&\dots &
\sigma^{n-1}(k_0)
\end{array}
\right]
\end{equation}
  By
taking $F$ to be various subfields of $\complexes$ containing
$\rationals(S)$ (the field generated
by the elements of $S$ over $\rationals$) in this formulation, and for each such $F$ taking
various $K$ and $\gamma$, a wide variety of space-time codes can be
constructed for a wide range of signal sets. 
For further simplicity of construction, particularly in the selection of the element $\gamma$ above, the authors of (\cite{SRS}) chose all their base fields $F$ to contain transcendental elements; in most cases, their cyclic extensions
$K/F$ were of the form $K_0(x)/F_0(x)$, where $K_0/F_0$ is a cyclic extension of number fields, and $x$ is a transcendental.  In these cases, the authors' construction yielded codes $X\colon S^{n^2} \rightarrow M_n(\complexes)$, i.e., with $k=n^2$.

Alamouti's
original code above arises as a special case of this formulation: the matrices of
Equation (\ref{eqn:alamouti_code}) are just the matrices of Equation (\ref{cycmat}) above specialized 
to the cyclic algebra
$(\complexes/\reals, \sigma, -1)$, where $\sigma$ stands for complex
conjugation.   This is nothing other than Hamilton's quaternions: the four-dimensional
$\reals$ algebra
$\reals \oplus \reals \imath \oplus \reals j \oplus \reals k$ subject to the relations 
$\imath^2 = j^2 =  -1$, $\imath j = -j \imath = k$. (The signal set in Alamouti's construction is
contained in $K$ instead of $F$, unless of course if $S$ is real.)

\section*{Satisfying the Coding Gain criterion} \label{secn:cgc}
The coding community immediately recognized the potential of
cyclic division algebras as a fundamental tool for constructing
space-time codes and began to work with the coding paradigm introduced in
\cite{SRS}. However, there was still a drawback.  While the specific codes of \cite{SRS} certainly
satisfied the rank criterion, their performance was not satisfactory.
The reason for this became clear: the
specific division algebras of \cite{SRS} were proposed only for
mathematical simplicity--merely as easy examples of the larger paradigm of division algebras--and were not optimized for the  coding
gain performance criterion above.  The use of
transcendental numbers in the codes in \cite{SRS} caused the determinants of the difference
matrices to
come arbitrarily close to zero and limited their performance.

This situation was quickly remedied in \cite{BRV} by a very clever
technique.  To provide a lower bound on the modulii of the
determinants of the difference of code matrices, the authors Belfiore, Rekaya
and Viterbo first constructed division algebras from cyclic
extensions $K/\rationals(\imath)$ 
and $\gamma \in \integers[\imath]$,
but then restricted the various
$k_i$ in the matrix (\ref{cycmat}) above to entries in $\mathcal{O}_K$, the ring of
integers of $K$.   The net result, as can easily be seen, is that
the determinant of the difference of any two such matrices  will live in
$\integers[\imath]$, and therefore will have modulus bounded below by $1$.
Moreover, this will be true no matter how large a subset of
$\integers[\imath]$ is used as the signal set.  They called this last 
property the ``nonvanishing determinant property'' and they called
the specific code they proposed the Golden Code.   It was so named
for the Golden Ratio that appears naturally: it is derived from
the division algebra
$(\rationals(\imath,\sqrt{5})/\rationals(\imath),\sigma, \imath)$.
Here, $\sigma$ is the automorphism of $K
=\rationals(\imath,\sqrt{5})$ that sends $\sqrt{5}$ to $-\sqrt{5}$
and acts as the identity on $\rationals(\imath)$. A
$\integers$-basis for $\mathcal{O}_K$ is given by $1$ and $\phi =
\frac{1+\sqrt{5}}{2}$. Write $\psi$ for $\sigma(\phi)
=\frac{1-\sqrt{5}}{2}$. For a signal set $S\subset
\integers[\imath] \subset \rationals(\imath) $ (the most common kind of signal set), this code
sends $S^4$ to $M_n(\complexes)$ via the matrix
\begin{equation}\label{GCMatrix}
\frac{1}{\sqrt{5}}
\left(%
\begin{array}{cc}
  s_{0,1}\alpha + s_{0,2}\alpha\phi& \imath(s_{1,1}\theta + s_{1,2}\theta\psi)\\
  s_{1,1}\alpha + s_{1,2}\alpha\phi & s_{0,1}\theta + s_{0,2}\theta\psi \\
\end{array}%
\right)
\end{equation}

Here, the $\dfrac{1}{\sqrt{5}}$ scale factor, $\alpha = {1+\imath(1-\phi)}$, and $\theta =
\sigma(\alpha) = {1+\imath(1-\psi)}$ are used to
shape the code (more on this ahead). Comparing with the matrix
(\ref{cycmat}) above and ignoring the scale factor, we see that $k_0 = s_{0,1}\alpha +
s_{0,2}\alpha\phi$ and $k_1 =s_{1,1}\alpha + s_{1,2}\alpha\phi$.
Note that this code encodes four information symbols in each matrix.  (A variant of this code, also based on the division
algebra $(\rationals(\imath,\sqrt{5})/\rationals(\imath),\sigma,
\imath)$, also incorporating the shaping criterion described ahead, 
is currently part of the IEEE 802.16e ``WIMAX''
standard.  The Alamouti code based on the quaternions is also part
of this standard.)

With the introduction of cyclic division algebras as a fundamental
construction paradigm and with the use of codes constructed with
entries from $\mathcal{O}_K$ for suitable extensions of
$\rationals(\imath)$, the subject of space-time coding took off. It is harmless and very
often actually useful to assume that the signal set $S$ is infinite:
typically, $S$ is assumed to be one of the standard lattices
$\integers$, $\integers[\imath]$ or the Eisenstein lattice
$\integers[\omega]$, where $\omega$ stands for the primitive third
root of unity $\frac{-1+\sqrt{-3}}{2}$. (Under these assumptions the
code forms an additive group, so one only needs to consider the rank
of $X(s_1,\dots,s_{k})$ and the modulus of the determinant $|\det
\theta X(s_1,\dots,s_{k})|$ in the rank criterion and the coding gain
criterion.)
 Coding theorists immediately
looked for specific constructions of division algebras of the form
$(K/F,\sigma,\gamma)$ for the cases where $F=\rationals$,
$F=\rationals(\imath)$, and $F=\rationals(\sqrt{-3})$,
corresponding to signal sets equaling one of the three lattices above.   While
such constructions have been known in principle to mathematicians
working with division algebras, the coding theorists absorbed the
necessary number-theoretic background in very short order and
explicitly constructed division algebras over such fields for all
indices $n$ (\cite{KR} and \cite{E5}).  (The hard task here is to
select $\gamma \in \mathcal{O}_F$ so that it has the property that $\gamma^i$ is not
a norm from $K$ to $F$ for $i=1,\dots,n-1$.) In all such cases, an
$\mathcal{O}_F$-basis $\beta_j$ of $\mathcal{O}_K$ is chosen, and each $k_i$ is
written as $\sum_{j=1}^n s_{i,j}\beta_j$ for $s_{i,j}$ in the
signal set. Thus, $n^2$ elements from the signal set are coded in
each matrix, and by construction, the determinant of each matrix
is nonzero and lies in one of the discrete lattices above. The
modulus of the determinant will therefore be bounded below by the
length of the shortest vector in the lattice so the code will
  have the nonvanishing determinant property.

\section*{Other Performance Measures} \label{secn:opm}
In parallel, as the subject became better understood, several
additional performance criteria started to be imposed on codes. In a
fundamental paper \cite{ZT}, Zheng and Tse provided a precise
quantification of the trade-off (known as the diversity-multiplexing
gain or ``DMG'' tradeoff) between information rate and reliability.   They
defined numerical measures for each of the benefits, and showed that
the pair of benefits lie in a region of the first quadrant whose
upper boundary is a piecewise linear concave up curve.  In the paper
\cite{E5} Vijay Kumar and his students showed that all codes
constructed from cyclic division algebras with the additional
nonvanishing determinant property will automatically perform at the
upper boundary of this region, and will hence be ``DMG optimal.''
This of course further cemented the use of cyclic division algebras
for code construction.

Another set of criteria were proposed by Oggier and coworkers in
the paper \cite{ORBE}.  One first rewrites the matrix
(\ref{cycmat})
 as a single $n^2\times 1$ vector.
When  $k_i =  \sum_{j=1}^n s_{i,j}\beta_j$ for $s_{i,j}$ in the
signal set and $\beta_j$ an $\mathcal{O}_F$ basis for
$\mathcal{O}_K$, this $n^2\times 1$ vector can be expressed as
   $M.\vbar$,
where $M$ is an $n^2 \times n^2$ matrix and $\vbar$ is the column
vector $(s_{0,1},s_{0,2},\dots, s_{i,j},\dots,s_{n-1,n})^T$.  One
now requires that the matrix $M$ be \textit{unitary} and that
$|\gamma|=1$.  The first condition is called ``good shaping'' and
the idea behind it is that this forces the average energy needed
to send the vector $\vbar$ without coding to be the same as that
needed to send it in the coded matrix form (\ref{cycmat}).  The
condition $|\gamma|=1$ causes the average energy transmitted per
antenna to be equal for all transmission. Oggier and coworkers
called such codes ``perfect'' and constructed perfect codes for
$n=2,3,4$ and $6$.  This was followed by work of Elia and
coworkers (\cite{ESV}) who constructed perfect codes for all
values of $n$, and additionally, showed that perfect codes satisfy
other information-theoretic properties such as
information-losslessness and approximate universality.

  The mathematics needed for the work on perfect
codes is quite interesting.  Analyzing the condition that $M$ be
unitary, we find that it is sufficient to make the following
matrix unitary:
$$
U(\{\beta_1,\dots,\beta_n\}) = \left[ \begin{array}{ccc} \beta_1 & \cdots & \beta_{n} \\
\sigma(\beta_1) & \cdots & \sigma(\beta_{n}) \\
 & \vdots & \\
\sigma^{n-1}(\beta_1) & \cdots  & \sigma^{n-1}(\beta_{n}) \\
\end{array} \right]
$$
Here, it is not necessary that the $\beta_j$ be an $\mathcal{O}_F$
basis of $\mathcal{O}_K$, it is sufficient that they be an
$\mathcal{O}_F$ linearly-independent subset of $\mathcal{O}_K$.
(So, for example, in the Golden Code (\ref{GCMatrix}) above,
$\alpha$ is chosen that with $\beta_1 = \alpha$ and $\beta_2 =
\alpha \phi$, the matrix $$
\left(%
\begin{array}{cc}
  \alpha & \alpha\phi \\
  \theta & \theta\psi \\
\end{array}%
\right)
$$ is unitary after being multiplied by the $\dfrac{1}{\sqrt{5}} $ scale factor.) So the question is: how to find  $\mathcal{O}_F$
submodules of $\mathcal{O}_K$ that satisfy this unitary condition?
 For $n=2^b$, it is easy to see that for the field $K =
 \rationals(\zeta)$ and $F = \rationals(\imath)$, where $\zeta$
 is a primitive $2^{b+2}$-th root of unity, the various powers of
 $\zeta$ are $\integers[\imath]$-linearly independent and satisfy the unitary condition above.
 For odd $n$, Elia and coworkers use a construction due to B. Erez
 (\cite{Erez})
that was needed in a different context: Erez was showing that for
certain cyclic extensions $K/\rationals$ with Galois group $G$,
the square-root of the inverse different is a free $\integers[G]$
module which has an orthogonal basis with respect to the usual
trace form on $K$ that sends $x,y$ to $Tr_{K/\rationals}(xy)$.

The most recent performance criteria for space-time codes, and in
some sense the most mathematically exciting, have come from Lahtonen and coworkers
 (\cite{HLRV}).  For the usual cases where $S$ is one of $\integers$,
$\integers[\imath]$,  $\integers[\omega]$, it is easy to see from the linearity of the code matrices $X$  that on  writing each $X$ as an $n^2\times 1$ vector as above and separating
the real and imaginary parts, one gets a full lattice in $\reals^{2n^2}$, i.e., the additive
group generated by $2n^2$ linearly independent vectors in $\reals^{2n^2}$. We refer to this lattice as the \textit{code lattice}.
  After normalizing
all code matrices  so that $\inf_{X\in\mathcal{X}}|\det(X)|=1$,
they
postulate that codes whose lattice points are the most dense in
$\reals^{2n^2}$ will have the best performance, and indeed, they find this
is borne out in several circumstances by simulations.
To obtain a suitable numerical measure for the relative density, they invert the situation: they normalize the code lattice to have fundamental volume $1$ instead.  Thus,
they
define the normalized minimum determinant of a code lattice
$\Lambda$ of rank $2n^2$ in a $\rationals(\imath)$ division
algebra of index $n$ (embedded in $M_n(\complexes)$) as the minimum of the modulii of the determinants $|\det(X(s_1,\dots, s_{n^2}))|$ as
$X(s_1,\dots, s_{n^2})$
runs through the lattice, divided by the fundamental
volume of $\Lambda$.
Since a smaller fundamental volume
represents a higher density,
the goal is to  construct codes whose code lattice $\Lambda$ would maximize this ratio among all full lattices in the division algebra.

Recall that if $D$ is a division algebra with center $F$ and if $R$ is a subring of $F$ whose quotient field is $F$, then an $R$-order in $D$ is a subring $T$ of $D$ containing $R$ that is finitely generated as an $R$-module and satisfies $TF = D$. A maximal $R$-order is one that is maximal with respect to inclusion.  In the typical situation where $S$ is one of  $\integers$,
$\integers[\imath]$, or $\integers[\omega]$, so $F$ is one of $\rationals$,
$\rationals(\imath)$, or $\rationals(\sqrt{-3})$, and where the $k_i$ of
 the matrices in (\ref{cycmat}) are
constrained to lie in $\mathcal{O}_K$ and $\gamma \in
\mathcal{O}_F$, the code matrices of (\ref{cycmat}) naturally form an $S$-order.  Thus the code matrices have a dual structure of  an $S$-order and a full lattice in $\reals^{2n^2}$.  Lahtonen and coworkers investigate the interplay between these two structures. They ask: how will the code's performance as measured by its normalized minimum determinant vary if, in addition to carrying its natural structure of a full $\integers$-lattice in $\reals^{2n^2}$, we choose our code to form an arbitrary $S$-order inside an $F$-division algebra?
In these cases, the minimum modulus of the determinants of the code matrices is
$1$, so it follows from the definition of the normalized minimum determinant that the smaller the fundamental volume of the lattice the better the code.  If $T_1$ and $ T_2$ are $S$-orders and $\Lambda_{T_1}$ and $\Lambda_{T_2}$ the corresponding lattices with fundamental volumes $V_{T_1}$ and $V_{T_2}$, then $T_1\subseteq T_2$ implies $\Lambda_{T_1} \subseteq \Lambda_{T_2}$, which in turn means that $V_{T_2} \le V_{T_1}$.  It follows therefore that the best normalized minimum determinant will arise when a maximal
order is used for the code. The authors then relate the fundamental volume of the code lattice
to the $\integers$-discriminant of the maximal order, and then invoke
known formulas for discriminants of maximal orders to compute the
best normalized minimum determinant of codes arising from
$\mathcal{O}_F$ orders inside a given division algebra.  In
particular, they show (for the fields $\rationals(\imath)$,
$\rationals(\sqrt{-3})$ and $\rationals$) that the best division
algebras to use will be ones that are ramified at precisely two of
the ``smallest'' primes of the field (where the size of a prime $P
= <\pi>$ is defined to be the modulus $|\pi|$).  Thus, for
$\rationals(\imath)$ for example, one needs to transmit on a code
arising from a maximal order inside a division algebra ramified
only at $(1+\imath)$ and $(2+\imath)$ (or $(2-\imath)$).  (Much of
this was part of Vehkalahti's Ph.D thesis.)

One of the drawbacks of using maximal orders is that the corresponding code lattice
may not have good shape. Thus, optimizing a code for minimum normalized determinant may destroy any optimization for shape.  
The recent
 work of Raj Kumar and
Caire  (\cite{CR}) proposes a very clever technique of mapping lattice points to certain
\textit{cosets} of a suitably chosen sublattice of a standard
cubic lattice; this smooths out
an irregular lattice and gives it better shape.   In particular, their technique applies to codes
from lattices from maximal orders and provides a further
performance boost in such cases.

\section*{Key Challenge: Decoding} \label{secn:dec}
What are some of the key problems that need to be solved in
space-time codes? Perhaps the biggest engineering challenge in the
subject is the issue of decoding. The problem quite simply is the
following: given the received vectors in $Y$ (see Equation
(\ref{eqn:block_model})), determine the entries of the matrix $X$
that represent the original information.
Assume that $k$ symbols are coded in the matrix $X$ and that the entries of $X$
are linear in the signal entries $s_1$, $\dots$, $s_k$ (typically arising from $\integers$, $\integers[\imath]$, or $\integers[\omega]$). 
By writing out sequentially the real and imaginary parts of each entry of $Y$, $W$, and 
$s_1$, $\dots$, $s_k$, we may rewrite Equation (\ref{eqn:block_model}) as
$\widetilde{Y} = Z\vbar + \widetilde{W}$.  Here $Z$ is an $2n^2\times 2k$ real
matrix that depends on $H$, $\theta$, and the parameters of the code matrix
$X$, $\vbar$ is the signal vector $(x_{1},y_{1}, \dots,
x_{i},y_{i},\dots,x_{k},y_{k})^T$ with $x_{i}$ and $y_{i}$ being the real and imaginary parts of $s_{i}$, and similarly for $\widetilde{Y}$ and $\widetilde{W}$.  If the columns of $Z$ were orthonormal, decoding would
be quite simple: we would have $Z^{T}\widetilde{Y} = \vbar
+Z^{T}\widetilde{W}$ with $Z^{T}\widetilde{W}$ also having independent, identically distributed Gaussian entries.  Hence, under maximum likelihood estimation,
$\vbar$ can be taken to be the closest vector  in $S^{k}$ (viewed inside the Euclidean space $\reals^{2k}$)
to $Z^{T}\widetilde{Y}$.  This is a very simple and computationally fast scheme: we march through  $Z^{T}\widetilde{Y}$ component pair by component pair and we find the element of the signal lattice $S$ closest to that component pair. 

The process above is called \textit{single symbol decoding.} (For $k < n^2$ this is the same as orthogonal projection on to the subspace of $\reals^{2n^2}$ determined by the columns of $Z$.)  
There are some nice situations where the matrix $Z$ is (essentially) orthogonal: this happens in the case of the Alamouti code, and more generally, in the codes satisfying Equation (\ref{eqn:orth_des}). The matrix $Z$ for such codes satisfies $Z Z^{T} = \theta^2 Tr(HH{^\H})I_{2k}$.  We may divide the relation $\widetilde{Y} = Z\vbar + \widetilde{W}$ by $\theta\sqrt{Tr(HH{^\H})}$.  The entries of the new noise vector $1/(\theta\sqrt{Tr(HH{^\H})})\widetilde{W}$ are still independent identically distributed Gaussian, while the columns of the matrix $1/(\theta\sqrt{Tr(HH{^\H})})Z$ are now orthonormal.   Thus single symbol decoding can be employed in all these cases.

But for other codes $Z$ is rarely orthogonal! In general, given that the entries of $W$ are independent identically distributed Gaussian, for maximum likelihood estimation one needs to search in $S^{n^2}$ (viewed inside the Euclidean space $\reals^{2n^2}$) for that vector $\vbar = (x_1,y_1,\dots, x_{n^2},y_{n^2})^T$ such that $Z\vbar$ is closest to $\widetilde{Y}$.  (Here we will
assume that $k=n^2$, as is usually the case for codes from cyclic division algebras.)  This can no longer be accomplished symbol by symbol, and one needs to search in the full space $S^{n^2}$ instead of just in $S$.  There is an algorithm called the \textit{sphere decoding algorithm} (see \cite{DGC} for instance) that accomplishes this search in an intelligent manner, but as is to be expected of any search in $S^{n^2}$, even this algorithm gets very cumbersome once $n$ exceeds $2$.  (However, in \cite{HV}, Hassibi and Vikalo show that under certain technical assumptions, the expected
complexity of the sphere decoding algorithm is polynomial, although the worst
case complexity is exponential.)

Since $Z$ is rarely orthogonal, we may ask  whether we can take advantage of the obvious
algebraic structure of the code and simplify the closest vector
problem for our particular application.
A very clever set of ideas of Lizzi et. al. (\cite{LRB}) does just that, and gives an approximate solution to
the decoding problem  for the Golden Code (Equation \ref{GCMatrix})
by reducing the situation to the action of $SL_2(\complexes)$ on
three dimensional hyperbolic space $\mathbb{H}^3$. Their work is a
veritable tour-de-force of the application of abstract mathematics
to engineering problems.  Their goal is to approximate the channel
matrix $H$ (normalized to have determinant $1$)
 by an element $U$ of determinant $1$ in the $\integers[\imath]$-order
$R= (\mathcal{O}_K/\integers[\imath], \sigma, \imath)$. Writing $H
= EU$ with $E$ simply being the error $HU^{-1}$, they argue that choosing $U$ so
that the Frobenius norm of $E^{-1}=UH^{-1}$ is minimized approximates the
original problem by the following: given a vector $Y$ in
$\complexes^{n^2}$ and an unknown vector $S$ in
$\integers[\imath]^{n^2}$ determine a ``best'' estimate of $S$ if
the difference vector $W = Y-S$ is known to be
\textit{approximately}  independent identically distributed Gaussian (in a suitable sense).  Given this assumption about the  the noise
vector $W$, a reasonable
way to proceed is to assume 
 that $W$
is actually independent identically distributed Gaussian.  In this
situation,  the maximum-likelihood
estimate of $S$ is obtained by taking the $i$-th entry of
$S$ to be the lattice point in $\integers[\imath]$ closest to the
$i$-th entry of $Y$.
  The authors find that their
scheme gives a fast and acceptably accurate decoding.

What is fascinating is the mathematics behind their choice of $U$.
First, they need to determine generators and relations for the
group of norm $1$ units $\U_1(R)$ of $R$ (i.e, the set of multiplicatively
invertible elements of $R$ whose
determinant as a code matrix is $1$). In general, it is very difficult to
find these for orders in division algebras, but in the case of
certain special quaternion algebras over number fields, generators
and relations for  $\U_1(R)$  is known.
Much of the ideas behind this goes
back to Poincare. The norm $1$ units in the order $R$ above
(modulo the subgroup $\{\pm 1\}$) turns out to be a Kleinian group,
i.e., a discrete subgroup of the projective special linear group
$PSL_2(\complexes)$.  As a subgroup of $PSL_2(\complexes)$,
$U_1(R)$ (modulo  $\{\pm 1\}$)
acts on the upper-half space
model of hyperbolic $3$-space $\mathbb{H}^3$  as a 
group of orientation-preserving isometries,
and Poincare's Fundamental Polyhedron Theorem gives a set of 
 generators and relations for such a group in terms of
certain automorphisms of a fundamental domain for the group.  Given a point $P$
in $\mathbb{H}^3$, the Dirichlet polyhedron centered on $P$ is the closure
of the set
of points $x$ such that $d_H(x,P) < d_H(g(x),P)$ for all $g\in U_1(R)$ (modulo $\{\pm 1\}$), $g\neq 1$,
where $d_H$ is the hyperbolic metric on
$\mathbb{H}^3$.  The authors construct a Dirichlet polyhedron centered on $J=(0,0,1)$; this is 
a fundamental domain for $U_1(R)$.  From this polyhedron, using Poincare's theorem and a computer search, 
they determine a set of
generators of $\U_1(R)$. They do this ahead of time, and store the
results.  Next, in real time, given a fading matrix $H$
(normalized to have determinant $1$), they need to find an element
$U$ of $\U_1(R)$ such that the Frobenius norm of $UH^{-1}$ is
minimized.  They observe that viewing $UH^{-1}$ as an element of
$PSL_2(\complexes)$ acting on $\mathbb{H}^3$, the Frobenius norm
of $UH^{-1}$ is just $2 \ cosh \  d_H(J,UH^{-1}(J))$, where $J$  and $d_H$
are as above.  Since $U$ is an isometry,
 they must find $U\in
\U_1(R)$ that minimizes $cosh\ d_H(U^{-1}(J), H^{-1}(J))$. From
the definition of Dirichlet polyhedra, this means that they need
to find a Dirichlet polyhedron centered on some $U^{-1}(J)$ which
contains $H^{-1}(J)$.  They use the geometry of $\mathbb{H}^3$ relative to
the action of $U_1(R)$ to find such a $U$: they just need to 
repeatedly consider the various 
 Dirichlet polyhedra
centered on $J$ and the various $g_i(J)$, where the $g_i$ run through the generators
of $\U_1(R)$ that they have computed ahead of time, along with
their inverses. 


\section*{Role of Mathematicians} \label{secn:role_math}
What is the role of mathematicians in this field?  The subject is
clearly very mathematical; yet, unlike classical coding theory
which now has a mathematical life of its own and can, for
instance, be thought of as a theory of subspaces of vector spaces
over finite fields, the center of gravity of space-time codes currently lies
very solidly in engineering. There is as yet no deep independent
``mathematics of space-time codes'': the driving force behind the
subject consists of fundamental engineering problems that need to
be solved before MIMO wireless communication reaches its full
practical potential, particularly for three or more antennas. 
This author therefore believes that, as things stand now, isolated
mathematical investigations of space-time codes that are not
grounded in concrete engineering questions would very likely lead
to sterile results.  At least for now, mathematicians can best
contribute to the subject by working in collaboration with
engineers who are
motivated by fundamental engineering questions.  This
author has found that the leading engineers in the field have a
practical and intuitive understanding of much abstract mathematics, but welcome help from trained mathematicians.  (This author has also found that they are
a genuine pleasure to collaborate with.)  There is clearly a lot of
work for mathematicians to do: particularly in decoding systems
with large numbers of receive and transmit antennas,  but also in other areas of MIMO communication 
that we have not touched upon in this article, such as cooperative
communication in networks, or noncoherent communication, where the matrix $H$ is not known to either the receiver or the transmitter.


\end{document}